\documentclass[11pt]{article}
\usepackage{amssymb}
=0cm \evensidemargin = 0cm 

\newtheorem{example}{Example}[section]
\newtheorem{defin}{Definition}[section]

\begin{document}
\title{\bf Pairing and Quantum Double of Finite Hopf C*-Algebras
\footnote{Supported by National Natural Science Foundation of China
(No.10301004)and Excellent Young Scholars Research Fund of Beijing
Institute of Technology(000Y07-25)}}
\author {\bf Ming LIU, Li Ning JIANG, Guo Sheng ZHANG
\\ Department of  Mathematics,\\
Beijing Institute of Technology, Beijing, 100081, P.R. China\\
E-mail: sunny$_{-}$sunabc@yahoo.com.cn}
\date{}
\maketitle

{\noindent \bf Abstract}  This paper defines a pairing of two
finite Hopf C*-algebras $A$ and $B$, and investigates the
interactions between them. If the pairing is non-degenerate, then
the quantum double construction is given. This construction yields
a new finite Hopf C*-algebra $D(A,B)$. The canonical embedding
maps of $A$ and $B$ into the double are both isometric.


{\noindent \bf MR (2000) Subject Classification} 46K70;16W30;81R05

\noindent Abbreviated title: Pairing of Hopf C*-algebras
\newpage

\section*{1\ \ \bf Introduction}

The Yang-Baxter equation firstly came up in a paper ([1]) as a
factorization condition of the scattering $S$ matrix in the
many-body problems in one dimension and in Baxter's work on
exactly solvable models in statistical mechanics. The equation
also plays an important role in the quantum inverses scattering
method created by Feddeev, Sklyamin and Takhtadjian for the
construction of quantum integrable systems. Since braided Hopf
algebras ([2]) can provide solutions for the Yang-Baxter equation,
attempts to find its solutions in a systematic way have led to the
construction of braided Hopf algebras, and moreover, led to the
theory of quantum group. Based on ([3]), Woronowicz ([4])
exhibited C*-algebra structures of quantum groups in the framework
of C*-algebra. From then on, the research on Hopf algebra has
always been going with that on C*-algebra. This leads to the
concept of Hopf C*-algebra ([5]). Indeed, by the Gelfand-Naimark
Theorem, an abelian C*-algebra can be understood as the space of
all complex continuous functions vanishing at infinity on a
locally compact space, and for this reason a C*-algebra can be
considered as a noncommutative locally compact quantum space.
Henceforth, a Hopf C*-algebra can be regarded as a noncommutative
locally compact quantum group.

In this paper, we are interested in finite Hopf C*-algebras.
Firstly, we propose the definition of Hopf *-algebra.

{\noindent\bf Definition 1.1} ([6, 7]) \textit{Let $A$ be a
*-algebra with a unit 1. Suppose that $\Delta:A\longrightarrow
A\otimes A$ is a *-homomorphism such that
$$(\Delta\otimes\iota)\Delta=(\iota\otimes\Delta)\Delta,$$
and $\varepsilon:A\longrightarrow\mathbb{C}$ is a *-homomorphism
such that
$$(\varepsilon\otimes\iota)\Delta=(\iota\otimes\varepsilon)\Delta=\iota,$$
where $\iota$ denotes the identity map. Finally, assume that
$S:A\longrightarrow A$ is a linear, anti-multiplicative map so
that for all $a\in A$, $$S(S(a)^{*})^{*}=a,$$
$$m(S\otimes\iota)\Delta(a)=m(\iota\otimes S)\Delta(a)=\varepsilon(a)1,$$
where $m:A\otimes A\longrightarrow A$ is the multiplication
defined by $m(a\otimes b)=ab$. Then $A$ is called a Hopf
*-algebra, and $\Delta, \varepsilon, S$ are called
comultiplication, counit and antipode respectively.}

Besides the assumption that $A$ is a finite dimensional Hopf
*-algebra, if $A$ is also a C*-algebra, $A$ is called a finite
Hopf C*-algebra, which satisfies the convolution property, i.e.,
$S^2=\iota$ naturally. And by ([8]), there exists an invariant
functional $\varphi$ on $A$ so that $\forall a\in A$,
$\varphi=\varphi\circ S$ and
$$(\varphi\otimes\iota)\Delta(a)=(\iota\otimes\varphi)\Delta(a)=\varphi(a)1.$$
After posing the notion of pairing of finite Hopf C*-algebras and
exhibiting the actions of dually paired Hopf C*-algebras on each
other, this paper gives the quantum double construction ([9]) out of
the underlying two Hopf C*-algebras and shows that the consequence
of this construction is again a finite Hopf C*-algebra with an
invariant integral.

Much of our work is inspired by the work of ([10, 11]). All algebras
in this paper will be algebras over the complex field $\mathbb{C}$.
Please refer to ([12]) for general results on Hopf algebra. In many
of our calculations, we use the standard Sweedler notation ([13]).
For instance, formula like
$m(S\otimes\iota)\Delta(a)=\varepsilon(a)1$ can be written as
$$\sum\limits_{(a)}S(a_{(1)})a_{(2)}=\varepsilon(a)1.$$

\section*{2\ \ \bf Pairing of Finite Hopf C*-Algebras}

In this section, we consider a bilinear form between finite Hopf
C*-algebras.

{\noindent\bf Definition 2.1} \textit{Let $A$ and $B$ be two
finite Hopf C*-algebras, and $<\cdot,\cdot>:A\otimes
B\longrightarrow \mathbb{C}$ be a bilinear form. Assume that they
satisfy: $\forall a_1, a_2, a\in A, b_1, b_2, b\in B$,
$$<\Delta(a), b_1\otimes b_2>=<a, b_1b_2>,$$
$$<a_1\otimes a_2, \Delta(b)>=<a_1a_2, b>,$$
$$<a^*, b>=\overline{<a, S_B(b)^*>},$$
$$<a, 1_B>=\varepsilon_A(a),$$
$$<1_A, b>=\varepsilon_B(b),$$
$$<S_A(a), b>=<a, S_B(b)>,$$
where $\varepsilon_A, S_A$ (resp. $\varepsilon_B, S_B$) denote the
counit and antipode on $A$ (resp. $B$) respectively. Then $(A, B,
<\cdot, \cdot>)$ is called a pairing of finite Hopf C*-algebras.}

{\noindent\bf Definition 2.2} ([14]) \textit{Suppose that $(A, B,
<\cdot, \cdot>)$ is a pairing of finite Hopf C*-algebras. If $B$
(resp. $A$) can separate the points of $A$ (resp. $B$) (i.e., if
$a_0\in A$ (resp. $b_0\in B$) such that $\forall b\in B$ (resp.
$a\in A$) $<a_0, b>=0$ (resp.$<a, b_0>=0$), then $a_0=0$ (resp.
$b_0=0$)) and we call the pairing is non-degenerate.}

Similar to the discussions in ([10]), we have the following results.
Firstly for any pairing of finite Hopf C*-algebras $(A, B, <\cdot,
\cdot>)$, we can define the linear mappings by using the standard
Sweedler notation:
$$\mu_{A,B}^{l}:A\otimes B\longrightarrow B, a\otimes b\mapsto\sum\limits_{(b)}b_{(1)}<a, b_{(2)}>,$$
$$\mu_{A,B}^{r}:B\otimes A\longrightarrow B, b\otimes a\mapsto\sum\limits_{(b)}<a, b_{(1)}>b_{(2)},$$
$$\mu_{B,A}^{l}:B\otimes A\longrightarrow A, b\otimes a\mapsto\sum\limits_{(a)}a_{(1)}<a_{(2)}, b>,$$
$$\mu_{B,A}^{r}:A\otimes B\longrightarrow B, a\otimes b\mapsto\sum\limits_{(a)}<a_{(1)}, b>a_{(2)}.$$
{\noindent\bf Proposition 2.3} \textit{Let $(A, B, <\cdot, \cdot>)$
be a pairing of finite Hopf C*-algebras. Then the maps
$\mu_{A,B}^{l}$ and $\mu_{A,B}^{r}$ are left and right actions of
$A$ on $B$, i.e., ($B$, $\mu_{A,B}^{l}$) is a left $A$-module and
($B$, $\mu_{A,B}^{r}$) is a right $A$-module, respectively.
Analogously, $\mu_{B,A}^{l}$ and $\mu_{B,A}^{r}$ are left and right
actions of $B$ on $A$, respectively.}

{\noindent \textit{Proof}} For all $a, a'\in A, b, b'\in B$, we
can obtain
$$
\begin {array}{llll}
\mu_{A,B}^{l}(aa'\otimes b)&=&\sum\limits_{(b)}b_{(1)}<aa', b_{(2)}>\\
&=&\sum\limits_{(b)}b_{(1)}<a, b_{(2)}><a', b_{(3)}>\\
&=&\sum\limits_{(b)}b_{(1)}<a, b_{(2)}<a', b_{(3)}>>\\
&=&\mu_{A,B}^{l}(a\otimes\mu_{A,B}^{l}(a'\otimes b)),
\end{array}
$$
which shows that ($B$, $\mu_{A,B}^{l}$) is a left $A$-module. In a
similar way, we can check other relations and we omit them here.
$\blacksquare$

For convenience, the previous actions will be denoted by
$\triangleright$ and $\triangleleft$ :
$$\mu_{A,B}^{l}(a\otimes b):=a\triangleright b,\ \ \mu_{A,B}^{r}(b\otimes a):=b\triangleleft a,$$
$$\mu_{B,A}^{l}(b\otimes a):=b\triangleright a,\ \ \mu_{B,A}^{r}(a\otimes b):=a\triangleleft b,$$
which mean ``$a$ acts from the left or right on $b$'' and ``$b$
acts from the left or right on $a$'' respectively, according to
the directions of the arrows $\triangleright$ and $\triangleleft$.

{\noindent\bf Lemma 2.4}  \textit{Let $(A, B, <\cdot, \cdot>)$ be
a pairing of finite Hopf C*-algebras. Then for all $a, a'\in A, b,
b'\in B$,
$$<b\triangleright a, b'>=<a, b'b>,\ \ \ <a\triangleleft b, b'>=<a, bb'>,$$
$$<a, a'\triangleright b>=<aa', b>,\ \ \ <a, b\triangleleft a'>=<a'a, b>.$$}
{\noindent \textit{Proof}}  From the implications of the notations
``$\triangleright$'' and ``$\triangleleft$'', the proof is
obvious. $\blacksquare$

From Lemma 2.4, we can get the following proposition at once.

{\noindent\bf Proposition 2.5}  \textit{Suppose that $(A, B, <\cdot,
\cdot>)$ is a non-degenerate paring of finite Hopf C*-algebras. Then
($A$, $\mu_{B,A}^{l}$, $\mu_{B,A}^{r}$) is a $B$-bimodule and ($B$,
$\mu_{A,B}^{l}$, $\mu_{A,B}^{r}$) is an $A$-bimodule.}

{\noindent \textit{Proof}}  Let $a\in A$ and $b_1, b_2, b_3\in B$,
and check $(b_1\triangleright a)\triangleleft b_2$ and
$b_1\triangleright (a\triangleleft b_2)$ paired with $b_3$. Using
Lemma 2.4, the associativity of $B$ implies
$$<(b_1\triangleright a)\triangleleft b_2, b_3>=<a, b_2b_3b_1>=
<b_1\triangleright (a\triangleleft b_2), b_3>,$$ which proves the
proposition for the non-degeneracy of the pairing.  $\blacksquare$

From Proposition 2.5, we will write $(b_1\triangleright
a)\triangleleft b_2$ as $b_1\triangleright a\triangleleft b_2$
briefly in sequence.

{\noindent\bf Remark 2.6} (1) The third axiom in Definition 2.1 is
also symmetric in $A$ and $B$:
$$\begin {array}{lll}
<a, b^*>&=&\overline{<a^*, S_B(b^*)^*>}\\
&=&\overline{<S_A^{-1}(a^*), S_B(S_B(b)^*)^*>}\\
&=&\overline{<S_A(a)^*, b>}.
\end{array}$$
(2) The last three axioms in Definition 2.1 are redundant if the
pairing $<\cdot, \cdot>$ is non-degenerate. From the first three
axioms of Definition 2.1 and Lemma 2.4, one can obtain for all $a,
a'\in A, b, b'\in B$, $<T_2^A(a\otimes a'), b\otimes b'>=<a\otimes
a',T_1^B(b\otimes b')>$, which also holds for the inverse mappings
of $T_2^A$ and $T_1^B$. Put $a=1_A$, $b=1_B$, by the
non-degeneracy of $<\cdot, \cdot>$, $<a'_{(1)}, 1_B>a'_{(2)}=a'$.
Applying $\varepsilon$ to the two sides of this equation yields
$<a', 1_B>=\varepsilon_A(a')$. Similarly, $<1_A,
b'>=\varepsilon_B(b')$. Using $<T_2^{A-1}(a\otimes a'), b\otimes
b'>=<a\otimes a',T_1^{B-1}(b\otimes b')>$ and Lemma 2.4, one can
get $<S_A(b'\triangleright a'), b\triangleleft a>=
<b'\triangleright a', S_B(b\triangleleft a)>,$ which implies the
last axiom.

\section*{3\ \ \bf The Quantum Double}

In what follows, we will only consider the action of $B$ on $A$,
where $A$ and $B$ are two dually paired finite Hopf C*-algebras. It
is easy to see that $A\otimes B$ can be made into a  linear space of
finite dimension in a natural way ([15]). Furthermore, we can turn
the linear space $A\otimes B$ into an associative algebra which has
an analogous algebra structure to the classical Drinfeld's quantum
double.

{\noindent \bf Definition 3.1} \textit{The quantum double $D(A,
B)$ of a non-degenerate paring of finite Hopf C*-algebras $(A, B,
<\cdot, \cdot>)$ is the algebra ($A\otimes B$, $m_D$) with the
multiplication map defined through
$$
\begin {array}{llll}
m_D((a, b)(a', b'))=\sum\limits_{(b)}(a
(b_{(3)})\triangleright a'\triangleleft (S^{-1}_{B}b_{(1)}), b_{(2)}b')\\
=\sum\limits_{(a')(b)}(aa'_{(2)}, b_{(2)}b') <a'_{(1)},
S^{-1}_{B}(b_{(3)})><a'_{(3)}, b_{(1)}>,
\end{array}
$$
where $(a, b), (a', b')$ are in the linear basis $B_D:=\{(a,b)\mid
a\in A, b\in B\}$ of $D(A, B)$.}

Following, we will write $m_D((a, b)(a', b'))$ as $(a, b)(a', b')$
directly. It is easy to see $(a, 1_B)(a', b)=(aa', b)$ and $(a,
b)(1_A, b')=(a, bb')$. In particular, $(1_A, 1_B)$ is the unit of
$D(A, B)$. Under the canonical embedding maps $i_A:a\mapsto (a,
1_B)$ and $i_B:b\mapsto (1_A, b)$, $A$ and $B$ become subalgebras
of $D(A, B)$.

{\noindent\bf Proposition 3.2}  \textit{The multiplication $m_D$
of the quantum double $D(A, B)$ is non-degenerate.}

{\noindent \textit{Proof}}  For a fixed element $(a, b)\in D(A,
B)$, suppose $(a, b)(a', b')=0$ for all $(a', b')\in D(A, B)$.
Particularly pick $a'=1_A$. Then $(a, b)(a', b')=(a, bb')=0$. If
$a\neq 0$, then $bb'=0$ for all $b'\in B$, which implies $b=0$ for
the non-degeneracy of the product on $B$. Thus we have $a=0$ or
$b=0$, i.e., $(a, b)=0$. Similarly one can prove that $(a, b)(a',
b')=0$ for all $(a, b)\in D(A, B)$ if and only if $(a', b')=0$.
$\blacksquare$

In order to avoid using too many brackets, we will use $Sa$ for
$S(a)$. On the basis $B_D$, set
$$*_D(a, b)=(a, b)^*:=\sum\limits_{(a)(b)}(a^*_{(2)}, b^*_{(2)})
<a^*_{(3)}, b_{(1)^*}><a^*_{(1)}, S^*_Bb_{(3)}>,$$ and extend it
anti-linearly to the whole space of $D(A, B)$. Then $(a,
1_B)^*=(a^*, 1_B)$, $(1_A, b)^*=(1_A, b^*)$. To describe the
*-structure of $D(A, B)$ exactly, we firstly do some preparing
work.

{\noindent\bf Lemma 3.3} \textit{$\forall (a, b)\in D(A, B)$, $(a,
b)^{**}=(a, b)$.}

{\noindent \textit{Proof}}
$$
\begin {array}{llll}
(a, b)^{**}\\=\sum\limits_{(a)(b)}(a^*_{(2)}, b^*_{(2)})^*<a_{(3)}, S_Bb_{(1)}><a_{(1)}, b_{(3)}>\\
=\sum\limits_{(a)(b)}(a_{(3)}, b_{(3)})\overline{<a^*_{(4)},
S^*_Bb_{(2)}><a^*_{(2)}, b^*_{(4)}>}
<a_{(5)}, S_Bb_{(1)}><a_{(1)}, b_{(5)}>\\
=\sum\limits_{(a)(b)}(a_{(3)}, b_{(3)})\overline{<a^*_{(4)},
S^*_Bb_{(2)}>}<a_{(2)}, S_Bb_{(4)}>
<a_{(5)}, S_Bb_{(1)}><a_{(1)}, b_{(5)}>\\
=\sum\limits_{(a)(b)}(a_{(3)}, b_{(3)})[<a_{(1)},
b_{(5)}><a_{(2)}, S_Bb_{(4)}>]
[<a_{(4)}, b_{(2)}><a_{(5)}, S_Bb_{(1)}>]\\
=\sum\limits_{(a)(b)}(a_{(2)}, b_{(3)})<a_{(1)},
b_{(5)}S_Bb_{(4)}>
<a_{(3)}, b_{(2)}S_Bb_{(1)}>\\
=\sum\limits_{(a)(b)}(a_{(2)},
b_{(2)})[\varepsilon_B(b_{(3)})<a_{(1)}, 1_B>]
[\varepsilon_B(b_{(1)})<a_{(3)}, 1_B>]\\
=\sum\limits_{(a)(b)}(a_{(2)}, b_{(1)})\varepsilon_B(b_{(2)})\varepsilon_A(a_{(1)})\varepsilon_A(a_{(3)})\\
=\sum\limits_{(a)}(a_{(2)}, b)\varepsilon_A(a_{(1)})\varepsilon_A(a_{(3)})\\
=(a, b),
\end{array}
$$
where we use relations $S_A((S_Aa)^*)^*=a$ and $<a^*,
b>=\overline{<a, S_Bb^*>}$ in the third and forth equations.
$\blacksquare$

{\noindent\bf Lemma 3.4} \textit{$\forall (a, b), (a', b')\in D(A,
B)$, $[(a, b)(a', b')]^{*}=(a', b')^*(a, b)^*$.}

{\noindent \textit{Proof}} We firstly prove the relation $[(1_A,
b)(a', b')]^{*}=(a', b')^*(1_A, b)^*$.

$$\begin {array}{llll}
[(1_A, b)(a', b')]^{*}\\
=[\sum\limits_{(a')(b)}(a'_{(2)}, b_{(2)}b')
<a'_{(1)}, S_Bb_{(3)}><a'_{(3)}, b_{(1)}>]^*\\
=\sum\limits_{(a')(b)}(a'_{(2)}, b_{(2)}b')^*\overline{<a'_{(1)}, S_Bb_{(3)}><a'_{(3)}, b_{(1)}>}\\
=\sum\limits_{(a')(b)(b')}(a'^*_{(3)},
b'^*_{(2)}b^*_{(3)})<a'^*_{(4)}, b'^*_{(1)}b^*_{(2)}>
<a'^*_{(2)}, S_Bb^*_{(4)}S_Bb'^*_{(3)}>\times\\ \overline{<a'_{(1)}, S_Bb_{(5)}><a'_{(5)}, b_{(1)}>}\\
=\sum\limits_{(a')(b)(b')}(a'^*_{(3)},
b'^*_{(2)}b^*_{(3)})\times\\ \overline{<a'_{(1)}, S_Bb_{(5)}>
<a'_{(2)}, b_{(4)}b'_{(3)}><a'_{(4)}, S_Bb'_{(1)}S_Bb_{(2)}><a'_{(5)}, b_{(1)}>}\\
=\sum\limits_{(a')(b)(b')}(a'^*_{(3)},
b'^*_{(2)}b^*_{(3)})\times\\ \overline{<a'_{(1)}, 1_B)>
<a'_{(2)}, b'_{(3)}><a'_{(4)}, S_Bb'_{(1)}S_Bb_{(2)}><a'_{(5)}, b_{(1)}>}\\
=\sum\limits_{(a')(b)(b')}(a'^*_{(3)},
b'^*_{(2)}b^*_{(3)})\times\\ \overline{<a'_{(1)}, 1_B)>
<a'_{(2)}, b'_{(3)}><a'_{(4)}, S_Bb'_{(1)}><a'_{(5)}, S_Bb_{(2)}><a'_{(6)}, b_{(1)}>}\\
=\sum\limits_{(a')(b')}(a'^*_{(3)},
b'^*_{(2)}b^*)\overline{<a'_{(1)}, 1_B>
<a'_{(2)}, b'_{(3)}><a'_{(4)}, S_Bb'_{(1)}><a'_{(5)}, 1_B>}\\
=\sum\limits_{(a')(b')}(a'^*_{(2)}, b'^*_{(2)})(1_A,
b)^*\overline{<a'_{(1)}, b'_{(3)}>
<a'_{(3)}, S_Bb'_{(1)}><a'_{(4)}, 1_B>}\\
=\sum\limits_{(a')(b')}(a'^*_{(2)}, b'^*_{(2)})(1_A,
b^*)\overline{<a'_{(1)}, b'_{(3)}>
<a'_{(3)}, S_Bb'_{(1)}>}\\
=(a', b')^*(1_A, b)^*.
\end{array}$$

Similarly, $(a, b)^*=[(a, 1_B)(1_A, b)]^*=(1_A, b)^*(a, 1_B)^*$
and then
$$\begin {array}{llll}
[(a, b)(a', b')]^{*}&=&[(a, 1_B)(1_A, b)(a', b')]^*\\
&=&[(1_A, b)(a', b')]^*(a, 1_B)^*\\
&=&(a', b')^*(1_A, b)^*(a, 1_B)^*\\
&=&(a', b')^*(a, b)^*,
\end{array}$$
which completes the proof. $\blacksquare$

Using Lemma 3.3 and Lemma 3.4, one can immediately get the
following result.

{\noindent\bf Proposition 3.5}  \textit{The involution $*_D$
renders $D(A, B)$ into a non-degenerate *-algebra.}

Furthermore, one can show that $D(A, B)$ has a Hopf *-algebra
structure. Indeed, under the following structure maps, $D(A, B)$
becomes a finite dimensional Hopf algebra naturally ([16]): $\forall
(a, b)\in D(A, B)$,
$$\Delta_D(a, b)=\sum\limits_{(a)(b)}(a_{(1)}, b_{(1)})\otimes(a_{(2)}, b_{(2)}),$$
$$\varepsilon_D(a, b)=\varepsilon_A(a)\varepsilon_B(b),$$
$$S_D(a, b)=\sum\limits_{(a)(b)}(S_Aa_{(2)}, S_Bb_{(2)})<a_{(1)}, S_Bb_{(3)}>
<a_{(3)}, b_{(1)}>.$$

{\noindent\bf Theorem 3.6} \textit{$D(A, B)$ is a Hopf *-algebra.}

{\noindent \textit{Proof}}  It suffices to show that $\Delta_D$
and $\varepsilon_D$ are *-homomorphisms and $\forall (a, b)\in
D(A, B)$, $S_D(S_D(a, b)^*)^*=(a, b)$.

(1) $\Delta_D$ is a *-homomorphism.
$$\begin {array}{llll}
\Delta_D((a, b)^*)&=&\Delta_D(((a, 1_B)(1_A, b))^*)\\
&=&\Delta_D((1_A, b^*)(a^*, 1_B))\\
&=&\Delta_D(1_A, b^*)\Delta_D(a^*, 1_B)\\
&=&\sum\limits_{(b)}(1_A, b^*_{(1)})\otimes(1_A, b^*_{(2)})
\sum\limits_{(a)}(a^*_{(1)}, 1_B)\otimes(a^*_{(2)}, 1_B)\\
&=&\sum\limits_{(a)(b)}(1_A, b_{(1)})^*(a_{(1)}, 1_B)^*\otimes(1_A, b_{(2)})^*(a_{(2)}, 1_B)^*\\
&=&\sum\limits_{(a)(b)}[(a_{(1)}, 1_B)(1_A, b_{(1)})]^*\otimes[(a_{(2)}, 1_B)(1_A, b_{(2)})]^*\\
&=&\sum\limits_{(a)(b)}(a_{(1)}, b_{(1)})^*\otimes(a_{(2)}, b_{(2)})^*\\
&=&\sum\limits_{(a)(b)}[(a_{(1)}, b_{(1)})\otimes(a_{(2)}, b_{(2)})]^*\\
&=&(\Delta_D(a, b))^*.
\end{array}$$
Similarly, $\varepsilon_D$ is a *-homomorphism.

(2) It is easy to see $S_D(a, 1_B)=(S_Aa, 1_B)$ and $S_D(1_A,
b)=(1_A, S_Bb)$. Thus
$$S_D(a, b)=S_D[(a, 1_B)(1_A, b)]=S_D(1_A,
b)S_D(a, 1_B)=(1_A, S_Bb)(S_Aa, 1_B),$$ and therefore,
$$(S_D(a,
b))^*=(S_Aa, 1_B)^*(1_A, S_Bb)^*=(S^*_Aa, S^*_Bb).$$ Using these
two relations, we have
$$\begin {array}{llll}
S_D(S_D(a, b)^*)\\=S_D(S^*_Aa, S^*_Bb)\\
=\sum\limits_{(S^*_Aa)(S^*_Bb)}(S_A(S^*_Aa)_{(2)},
S_B(S^*_Bb)_{(2)})
\times\\<S^*_Aa_{(1)}, S_B(S^*_Bb)_{(3)}><S^*_Aa_{(3)}, S^*_Bb_{(1)}>\\
=\sum\limits_{(a)(b)}(a^*_{(2)}, b^*_{(2)})<a^*_{(3)}, b^*_{(1)}><S^*_Aa_{(1)}, S_B(S^*_Bb)_{(3)}>\\
=\sum\limits_{(a)(b)}(a^*_{(2)}, b^*_{(2)})<a^*_{(3)}, b^*_{(1)}><a^*_{(1)}, S_Bb^*_{(3)}>\\
=(a, b)^*.
\end{array}$$
$\blacksquare$

Now it is time to consider the C*-algebra structure of $D(A, B)$.

{\noindent\bf Lemma 3.7} \textit{Let $\varphi_A$ and $\varphi_B$
be invariant integrals on $A$ and $B$, respectively. $\forall (a,
b)\in D(A, B)$, set
$$\theta((a, b)):=\varphi_A(a)\varphi_B(b).$$
Then $\theta$ is a faithful positive linear functional on $D(A,
B)$. }

{\noindent \textit{Proof}}  $(a, b)(a, b)^*=(a, b)(1_A, b^*)(a^*,
1_B) =(a, bb^*)(a^*, 1_B)$. In the following, we denote $bb^*$ by
$c$ briefly.
$$\begin {array}{llll}
\theta((a, b)(a, b)^*)&=&\theta((a, c)(a^*, 1_B))\\
&=&\sum\limits_{(a)(c)}\theta((aa^*_{(2)}, c_{(2)}))<a^*_{(1)}, S_Bc_{(3)}><a^*_{(3)}, c_{(1)}>\\
&=&\sum\limits_{(a)(c)}\varphi_A(aa^*_{(2)})\varphi_B(c_{(2)})<a^*_{(1)}, S_Bc_{(3)}><a^*_{(3)}, c_{(1)}>\\
&=&\sum\limits_{(a)(c)}\varphi_A(aa^*_{(2)})<a^*_{(1)}, S_Bc_{(3)}><a^*_{(3)}, \varphi_B(c_{(2)})c_{(1)}>\\
&=&\sum\limits_{(a)(c)}\varphi_A(aa^*_{(2)})<a^*_{(1)}, S_Bc_{(2)}><a^*_{(3)}, 1_B>\varphi_B(c_{(1)})\\
&=&\sum\limits_{(a)(c)}\varphi_A(a\varepsilon_A(a^*_{(3)})a^*_{(2)})
<a^*_{(1)}, S_Bc_{(2)}>\varphi_B(c_{(1)})\\
&=&\sum\limits_{(a)(c)}\varphi_A(aa^*_{(2)})<a^*_{(1)}, S_Bc_{(2)}>\varphi_B(c_{(1)})\\
&=&\sum\limits_{(a)(c)}\varphi_A(aa^*_{(2)})<a^*_{(1)}, \varphi_B\circ S_Bc_{(1)}S_Bc_{(2)}>\\
&=&\sum\limits_{(a)}\varphi_A(aa^*_{(2)})<a^*_{(1)}, 1_B>\varphi_B(c)\\
&=&\sum\limits_{(a)}\varphi_A(aa^*_{(2)})\varepsilon_A(a^*_{(1)})\varphi_B(c)\\
&=&\varphi_A(aa^*)\varphi_B(c)\geq 0,
\end{array}$$
where we use the relation $\varphi_B\circ S_B=\varphi_B$ for the
last third and forth equations.

 It is clear that $\theta((a, b)(a, b)^*)=0$ if and only if
$a=0$ or $b=0$, which implies $(a, b)=0$. Thus $\theta$ is a
faithful positive linear functional on $D(A, B)$. $\blacksquare$

{\noindent\bf Theorem 3.8} \textit{$D(A, B)$ is a finite Hopf
C*-algebra.}

{\noindent \textit{Proof}}  Using the result in Lemma 3.7, one can
construct the associated GNS representation of $D(A, B)$ ([17]):
$\forall x,y \in D(A, B)$, set
$$<x, y>_{_{\theta}}=\theta(y^*x),$$
where $<x, y>_{_{\theta}}$ denotes the inner product of $x$ and
$y$. Thus $D(A, B)$ turns into a Hilbert space $K$. For $d\in D(A,
B)$, define
$$\pi(d): K\longrightarrow K,\ \  x\mapsto dx.$$
Using ([17]), $(\pi, K)$ is a faithful *-representation of $D(A,
B)$, and hence $D(A, B)$ can embeds into $B(K)$ isometrically
through
$$\pi:D(A, B)\longrightarrow B(K),\ \  d\mapsto \pi(d).$$
Again $D(A, B)$ is finite dimensional, therefore, it is a
C*-algebra with C*-norm $\|(a, b)\|=(\theta((a, b)(a,
b)^*))^{1/2}$. $\blacksquare$

{\noindent\bf Remark 3.9} A short calculation shows that $\theta$
coincides with $\varphi_A\otimes\varphi_B$, which is indeed an
integral on $D(A, B)$. Using the relation $\theta((a, b)(a,
b)^*)=\varphi_A(aa^*)\varphi_B(bb^*)$, one can get $\|(a,
b)\|=\|a\|\|b\|$. In particular, $\|(a, 1_B)\|=\|a\|$ (resp.
$\|(1_A, b)\|=\|b\|$), which implies that the canonical embedding
map $i_A$ (resp. $i_B$) is isometric.

{\noindent\bf Example 3.10} Let $H$ be a finite Hopf C*-algebra and
$H'$ be its dual, which is also a finite Hopf C*-algebra by ([8]).
They are naturally dually pairing and have invariant integrals,
denoted by $h$ and $h'$ respectively. Drinfeld's quantum double
$D(H)$ of $H$, which is defined as the bicrossed product of $H$ and
$H'$, is a special case of our construction. One ([11]) can prove
that it is also a finite Hopf C*-algebra and has an invariant
integral $h\otimes h'$.

\end{document}